\pdfoutput=1
\RequirePackage{ifpdf}
\ifpdf 
\documentclass[pdftex]{sigma}
\else
\documentclass{sigma}
\fi

\numberwithin{equation}{section}

\newtheorem*{Theorem*}{Theorem}

 { \theoremstyle{definition}

 }

\begin{document}
\allowdisplaybreaks

\renewcommand{\thefootnote}{}

\renewcommand{\PaperNumber}{092}

\FirstPageHeading

\ShortArticleName{Seven Concepts Attributed to Sim\'eon-Denis Poisson}

\ArticleName{Seven Concepts Attributed to Sim\'eon-Denis Poisson\footnote{This paper is a~contribution to the Special Issue on Differential Geometry Inspired by Mathematical Physics in honor of Jean-Pierre Bourguignon for his 75th birthday. The~full collection is available at \href{https://www.emis.de/journals/SIGMA/Bourguignon.html}{https://www.emis.de/journals/SIGMA/Bourguignon.html}}}

\Author{Yvette KOSMANN-SCHWARZBACH}

\AuthorNameForHeading{Y.~Kosmann-Schwarzbach}

\Address{Paris, France}
\Email{\href{mailto:yks@math.cnrs.fr}{yks@math.cnrs.fr}}

\ArticleDates{Received September 21, 2022, in final form November 25, 2022; Published online November 29, 2022}

\Abstract{Sim\'eon-Denis Poisson was 25 years old when he was appointed Professor of Mathematics at the \'Ecole Polytechnique in 1806. Elected to the Paris Acad\'emie des Sciences six years later, he soon became one of its most influential members. The origin and later developments of the many concepts in mathematics and physics that bear his name make interesting stories, a few of which we shall attempt to sketch in this paper.}


\Classification{01A55; 01A70; 31A30; 31J05; 60G55}

\begin{flushright}
\begin{minipage}{59mm}
\it Pour Jean-Pierre Bourguignon\\ \`a l'occasion de son 75$^{\,e}$ anniversaire
\end{minipage}
\end{flushright}

\renewcommand{\thefootnote}{\arabic{footnote}}
\setcounter{footnote}{0}

\section{The mathematician and physicist Poisson (1781--1840)}

Pithiviers is a small town in France, 50 miles south of Paris, renowned for
a special, tasty kind of pastry, called a ``pithiviers'', and for the high-quality
honey from the neighboring countryside. But it has another claim to fame.
It was the birthplace in 1781 of Sim\'eon-Denis Poisson, who would become
the mathematician and mathematical physicist whose name is attached to
the Poisson distribution, the Poisson brackets, Poisson geometry, Poisson
algebras and many other concepts, formulas, equations and theorems.

He was just eight years old when the French Revolution broke out. His father was
a retired soldier who held a modest administrative position. The Revolution
allowed boys from families such as his to obtain a decent education. In 1794,
the \'Ecole Polytechnique, first called \'Ecole centrale des travaux publics, was
created to prepare engineers with a solid scientific background. Then, as
now, admission was by a competitive examination. Encouraged by one of
his teachers in high school, and equipped with a certificate attesting to his
deep love of liberty, equality and all the fundamental beliefs of the Republic,
including a ``hate for tyrants'', Poisson sat for the entrance examination
in 1798 and took first place. This was the beginning of a very brilliant
career, through many regime changes, first the revolutionary Republic, then
Napoleon's Empire, the Restauration of the royalty in 1815, under Louis
XVIII until his death in 1824, followed by the more autocratic Charles~X, then
the revolution of 1830 and the constitutional monarchy of Louis-Philippe.
Poisson did not live to see France's later regime, the short-lived republic of
1848, since he died in 1840 at the age of 58.

Ten years after his death, a lifesize statue of Poisson was installed in his
hometown. To this day, there is a square in the center of Pithiviers that bears
the name ``Place Denis Poisson'', but the brass statue has disappeared, like
so many others in France and elsewhere, having been melted down during
the German occupation of Pithiviers in the Second World War.

When Poisson entered the \'Ecole Polytechnique, the professors were the
most distinguished scientists of the time, Joseph Louis Lagrange (1736--1813)
and Gaspard Monge (1746--1818) taught mathematics, Pierre-Simon Laplace
(1749--1827) was an examiner for mathematics, Jean Baptiste Biot (1774--1862) was an examiner for physics, and Antoine-François de Fourcroy (1755--1809) was professor of chemistry.

Volume 4 of the {\it Journal de l'\'Ecole polytechnique}, dated 1801--1802, contains
three papers written by Poisson in 1799 or early 1800, when he was still
a student. One deals with the classification of quadrics and is an “Addition”
to an article by Monge and Jean Nicolas Pierre Hachette (1769--1834) who
were both professors at the \'Ecole. This 3-page note \cite{hachette} must have been a
remark that he made on a publication of his teachers, and it is signed by one
of them, Hachette, together with Poisson. It is both the first and last paper
ever co-authored by Poisson.

In the same volume, we find the paper ``Essay on elimination in algebraic
equations'',\footnote{M\'emoire sur l'\'elimination dans les \'equations alg\'ebriques~\cite{equations}.
} an essay on the elimination of variables in systems of algebraic
equations, which contains a new, simplified proof of the theorem of \'Etienne
B\'ezout (1730--1783) on the degree of the resultant attached to a pair of algebraic
curves, and the ``Essay on the plurality of integrals in the calculus
of differences'',\footnote{M\'emoire sur la pluralit\'e des int\'egrales dans le calcul des
diff\'erences~\cite{calcul}.} that he had read before the Institut National, which had replaced
the Acad\'emie des Sciences in 1795, on the 16th of Frimaire in the year
9 of the French revolutionary calendar, i.e., December 8, 1800. In that paper,
he generalized a remark of Laplace on the solutions of first-order difference
equations. The report by two members of the Academy, Adrien-Marie Legendre
(1752--1833) and Sylvestre François Lacroix (1765--1843), is preserved
in the Archives de l'Acad\'emie des Sciences. In the conclusion of their four-page
long ``Report on a paper by Citizen Poisson on the number of complete
integrals of which equations of finite differences are susceptible'',\footnote{Rapport
sur un M\'emoire du Citoyen Poisson relatif au nombre d'int\'egrales complettes
[sic] dont les \'equations aux diff\'erences finies sont susceptibles.} they wrote:
``It follows at a minimum that the theory established by this young geometer
is correct, and even though it may not be susceptible to useful applications
in the problems that lead to this type of equations, one must always regard
the clarification of a problem of analysis which, until the present, remained
in great obscurity as contributing to the progress of science'',\footnote{Il r\'esulte
au moins que la th\'eorie \'etablie par ce jeune g\'eom\`etre est exacte, et quand
m\^eme elle ne serait pas susceptible d'applications utiles dans les probl\`emes qui conduisent
\`a ce genre d'\'equations, on doit toujours regarder commme contribuant au progr\`es de la
science, l'\'eclairicssement d'un probl\`eme d'analyse qui jusqu'\`a pr\'esent \'etait rest\'e dans une
grande obscurit\'e.} and they recommended
the publication of Poisson's paper. Legendre and Lacroix were
not enthusiastic but they certainly did not discourage the promising young
mathematician.

Upon finishing his studies at the \'Ecole, Poisson was immediately appointed
as an assistant, and in 1802 he was called upon to ``take over temporarily
the duties of Citizen Fourier'', Joseph Fourier (1768--1830) who would soon develop
the Fourier series and integrals. Four years later, at the age of 25, he was
appointed ``Instituteur d'Analyse'', i.e., full professor of mathematics, to replace
Fourier for whom he had already been substituting for four years. The
archives of the \'Ecole Polytechnique contain the confirmation of Poisson's
nomination to replace Fourier, dated 11 Brumaire, year 11 (November~2,
1802), as well as the covering letter, dated March~17, 1806, of the Emperor's
official decree naming Poisson ``Instituteur d'Analyse'' at the \'Ecole Polytechnique.

In 1804, Poisson appears as a handsome young professor in a portrait by
the painter E.~Marcellot, that is now in the collection of \'Ecole Polytechnique.

In 1809, Napoleon decreed the opening of a re-organized ``Universit\'e
Imp\'eriale'', and Poisson was named its first professor of mechanics. A poster
announcing the opening of classes in April 1811 may still be seen in the
collection of the Biblioth\`eque Nationale de France and it specifies that Poisson's
lectures at the Sorbonne, to use the old name that has survived all
reforms, student revolts and re-organizations to this day, would be delivered
on Mondays and Fridays. (The other professors were Lacroix, for differential
and integral calculus, Louis-Benjamin Franc\oe ur (1773--1849) for advanced
algebra, Hachette for descriptive geometry, and Biot for astronomy.)

When there was finally a vacancy in the Academy of Sciences it was in
the physics section in 1812 and Poisson was elected to fill it. His role in the
Academy was soon preeminent, particularly in the election of new members.
He was charged with writing numerous reports. One, written in 1816, I~find
particularly interesting because it shows Poisson, already a respected member
of the Academy, in a position to judge an unknown, young scientist, just as
he had been judged by Legendre and Lacroix, back in 1800. In it, together
with Andr\'e-Marie Amp\`ere (1775--1836), he reports on a paper by Claude
Pouillet (1790--1868) on the phenomenon of colored rings. Their conclusion
on the work submitted by this ``young physicist'' is very reminiscent of that of
Legendre and Lacroix on the ``young geometer'' Poisson, sixteen years earlier,
and they too recommend the publication of the essay. Their judgment was
fair, since Pouillet went on to teach at Polytechnique and at the Sorbonne,
and to be elected to the Academy. This report in the Archives de l'Acad\'emie
des Sciences is an autograph manuscript of Poisson.

Another aspect of Poisson's role in the Acad\'emie was important even
if it was sometimes controversial. He often had to read papers submitted
for publication to the \textit{M\'emoires de l'Acad\'emie}. In several instances he rendered
a great service to the scientific community by summarizing in his own
terms the main points of these papers, before, sometimes many years before,
they were revised and published. A controversy arose when Poisson made
use of such material, with insufficient acknowledgement of his source. The
best known case is his bitter dispute with Fourier. Having had access to the
manuscript on the theory of heat that Fourier had submitted to the Academy
in 1807, he published a detailed account of the subject in 1808 in the {\it Bulletin
de la Soci\'et\'e Philomatique} ({\it Bulletin des Sciences}), under the title “On the
propagation of heat in solid bodies, by M.~Fourier''\footnote{M\'emoire sur la
propagation de la chaleur dans les corps solides, par M. Fourier.} and signed~P, his initial.
While, from 1814 to 1825, Poisson published many essays on trigonometric
series and the theory of heat, Fourier's text, revised, would be published in
1822~\cite{fourier}, much later than Poisson's first papers. A virulent debate over precedence
broke out in 1815 between the two scientists, during which Fourier
wrote a letter to Laplace, blaming both Poisson and Biot: ``[They] recognize
that they could not obtain up to now any results different from mine [\dots] but
they say that they have another method for formulating them, and that this
method is excellent and the true one. [\dots] But it does not extend the limits
of science to present results that one has not found in a form that one says
is different''.\footnote{[Ils] reconnaissent qu'ils n'ont pu donner jusqu'ici aucun r\'esultat diff\'erent des miens
[\dots] mais ils disent qu'ils ont une autre mani\`ere de les exposer et que cette mani\`ere est
excellente et la v\'eritable. [\dots] Mais ce n'est pas reculer les limites des sciences que de
pr\'esenter sous une forme que l'on dit \^etre diff\'erente des r\'esultats que l'on n'a pas trouv\'es
soi-m\^eme.}

When Gaston Darboux (1842--1917) edited the works of Fourier in 1890,
he included Poisson's 1808 account, explaining, ``This article [\dots] is not by
Fourier. Signed with the initial ``P'', it was written by Poisson who was an
editor of the mathematical portion of the {\it Bulletin des Sciences}. Because
of the historic interest that it presents as the first publication that made
Fourier's theory known, we believed that we should reproduce it in its entirety''.\footnote{Cet
Article [\dots] n'est pas de Fourier. Sign\'e de l'initiale P, il a \'et\'e r\'edig\'e par Poisson
qui \'etait un r\'edacteur du {\it Bulletin des Sciences} pour la partie math\'ematique. \`A raison de
l'int\'er\^et historique qu'il pr\'esente comme \'etant le premier \'ecrit où l'on ait fait connaître la
th\'eorie de Fourier, nous avons cru devoir le reproduire int\'egralement.}
In 1808, Fourier had derived the heat equation and solved it in a
particular case, by expanding the sought-after solution into a trigonometric
series. Poisson established the heat equation for the case of variable conductivity
in his ``Essay on the distribution of heat in solid bodies'',\footnote{M\'emoire sur la distribution
de la chaleur dans les corps solides.} published
in the {\it Journal de l'\'Ecole polytechnique} in 1823, and included it in his 1835
book, ``Mathematical Theory of Heat''\footnote{Th\'eorie math\'ematique de la chaleur.}~\cite{21}.

From 1820 until his death, Poisson played an important role in the organization
of education in France, as a member and then, after 1822, as
treasurer of the Royal Council of Public Education.\footnote{Conseil royal de l'instruction publique.} As head of mathematics
in France, he wielded considerable influence and used his authority to do
battle to insure the teaching of mathematics to all students, including those
primarily studying humanities. As president of the jury of the ``agr\'egation''
for the selection of teachers, he strove to maintain a single examination for
both mathematics and physics. He championed mathematics but he also
understood that the two fields had to be developed simultaneously at the
middle and high school levels and at the university, just as they were in his research.

\section{Some stories}

Poisson's name appears in so many contexts in mathematics and in physics
that discovering the earliest formulation of each of his concepts or theorems
in his more that 200 publications would be a very time-consuming enterprise.
Here, I shall mention only a~few of them.

\subsection{Poisson's equation in electromagnetism}

Poisson's first important contributions to the theory of electrostatics date
from 1811 to 1813, when he took up the problem of determining the distribution
of electrical charges in charged bodies by applying analytical techniques,
such as series expansions, while he began treating magnetism in an article
published later. In his ``Essay on the theory of magnetism in motion''\footnote{M\'emoire sur
la th\'eorie du magn\'etisme en mouvement.}~\cite{18},
which appeared in the {\it M\'emoires de l'Acad\'emie des Sciences} of 1823, one
finds, on p. 463, ``$\Delta V=0$, $=-2k\pi$, $=-4k\pi$, depending on whether point
$M$ is located outside, on the surface of, or inside the volume in question''.\footnote{$\Delta V=0$, $=-2k\pi$, $=-4k\pi$,
selon que le point M sera situ\'e en dehors, \`a la surface ou
en dedans du volume que l'on consid\`ere”.} Here $k$ denotes the constant charge density of the charged body. (It was
Friedrich Gauss (1777--1855) who would later treat the case of variable density.)
The case of the equation $\Delta V=0$ was already well known, being the
famous Laplace equation. Poisson had derived an equation satisfied by the
potential at a point interior to the charged body, but the novelty in the 1823
paper was to treat the case of a point on the surface of the body. This case
became known as ``Poisson's equation'' in electromagnetism. It was indeed
discovered by Poisson.

The importance of these articles was immediately recognized by George
Green (1793--1841), after whom are named the Green function and the
Green--Riemann theorem. In 1828, in the preface to his {\it Essay on the Application of
Mathematical Analysis to the Theories of Electricity and Magnetism,} Green
cited Poisson's work prominently and, speaking of the papers of 1811 and
1812, he wrote, ``Little appears to have been effected in the mathematical
theory of electricity [\dots] when M. Poisson presented to the French Institut
two memoirs of singular elegance, relative to the distribution of electricity on
the surfaces of conducting spheres, previously electrified and put in presence
of each other''.

The fame of Poisson's mathematical theory of electrostatics is reflected
in the judgment of E.T.~Whittaker (1873--1956) in his {\it History of the Theories
of Æther and Electricity} (1910). Regarding Poisson's 1812 essay he
wrote, ``Electrostatical theory was, however, suddenly advanced to quite a
mature state of development by Sim\'eon-Denis Poisson, in a memoir which
was read to the French academy in 1812 [\dots]. The rapidity with which in a
single memoir Poisson passed from the barest elements of the subject to such
recondite problems as those just mentioned may well excite admiration''. He
concluded: ``His success is, no doubt, partly explained by the high state
of development to which analysis had been advanced by the great mathematicians
of the eighteenth century; but [\dots] Poisson's investigation must be
accounted a splendid memorial of his genius''. Later, he examined ``Poisson's
theory of magnetic induction'', rejecting his interpretation of the physics of
the situation but noting that the formulas derived by Poisson are valid.

And later, in 1939, the historian of science, Mario Gliozzi, in a paper analyzing
``Poisson's contribution to electricity theory'',\footnote{Il contributo del Poisson all'electrologia.}
concluded that Poisson's 1813 publication was a most remarkable paper.

\subsection{Poisson's ratio}

The story of Poisson's ratio is that of a concept whose present, everyday
applications are as surprising as they are numerous. I heard a beautiful talk
by Tadashi Tokieda in 2012 in Paris, that started with Poisson, continued
with origami, and went on to an astonishing variety of contemporary questions
in materials science, and I later read the article by the physicist George
Neville Greaves (1945--2019), ``Poisson's ratio over two centuries: challenging
hypotheses'' \cite{greaves}, which gave both a historical account and a detailed description
of the current theory and practice of this concept, which I shall briefly
summarize here. It started with a ``shape versus volume concept'', a hint already
given by Poisson in his early {\it Trait\'e de M\'ecanique}~\cite{14}, first published
in 1811, where he wrote on page~176 of volume~2:
\begin{quotation}
\noindent For each of the elements into which we have divided the amount of fluid matter, its shape
will be altered during the time ${\rm d}t$, and also its volume will change if the fluid is compressible;
but, since its mass must remain unaltered, it follows that, if we seek to determine its volume and its
density at the end of time $t + {\rm d}t$, their product will necessarily be the same as after time~$t$.\footnote{Chacun des \'el\'ements dans lesquels nous avons partag\'e la masse fluide, changera
de forme pendant l'instant~${\rm d}t$, et il changera m\^eme de volume, si le fluide est compressible;
mais comme sa masse devra toujours rester la m\^eme, il s'ensuit que, si nous cherchons ce que
deviennent son volume et sa densit\'e \`a la fin du temps $t + {\rm d}t$, leur produit devra \^etre
le m\^eme qu'\`a la fin du temps $t$.}
\end{quotation}

In a Note in the \textit{Annales de chimie et de physique} in 1827, ``On the extension of elastic threads and
plates''\footnote{Sur l'extension des fils et des plaques \'elastiques.}~\cite{19}, Poisson introduced
the dimensionless ratio that bears his name, and, by means of a computation based on Laplace's theory of
molecular interaction, announced that its value is~$\frac{1}{2}$, in accordance with recent experiments on
brass [in French, {\it{laiton}}] by baron Charles Cagniard de La Tour (1777--1859) and F\'elix Savart (1791--1841),
on the vibrations of plates, whose results had been recently presented to the Academy. Poisson further developped
the theory of elasticity in several papers that he read before the Acad\'emie des Sciences between 1823 and 1828
and published from 1828 to 1830, introducing the ratio of the strain in the transverse direction to the strain
in the primary direction.

About ten years later, as the precision of the measurements in experiments increased, the constancy of Poisson's
ratio for all materials was proved to be wrong, but the conflict between the hypothesis of interacting molecules and the continuum theory
of Sophie Germain (1776--1831) and Augustin Cauchy (1789--1857) was not resolved until much later. In the 1860's, James Clerk Maxwell
(1831--1879) defended Poisson's viewpoint, but William Thomson (Lord Kelvin, 1824--1907) declared that it had already been proved false
by George Stokes (1819--1903) in 1845.
Until the 1970's the variability of Poisson's ratio with every kind of material had only been established experimentally
by engineers ``for whom macroscopic properties were sovereign''~\cite{greaves}. Still, it was shown that its variability, unlike
that of the other elastic moduli, is restricted within the range [$-1,\frac{1}{2}$].
The number of publications concerning the Poisson ratio increased exponentially after 1970, when it was discovered that
this concept helped to understand ``the narrowing of arteries during hypertension, the resilience of bones and medical
implants, the rheology of liquid crystals, the shaping of ocean floors, the oblateness of the Earth, and planetary
seismology after meteor impact''~\cite{greaves}.
This is when materials with negative Poisson's ratio began to appear and found countless, important applications that were aptly
evoked in Tokieda's entertaining lecture. This was the result of the work of many physicists, beginning with Roderic Lakes in
1987 and including Greaves, the author of the article that we have attempted to summarize here, who observed that
``Sim\'eon-Denis Poisson is particularly remembered for a ratio, a dimensionless quantity, which today has surprisingly acquired
a ubiquitous physical significance''.\footnote{See also the comprehensive article by G.N.~Greaves et al., ``Poisson's ratio and
modern materials'', {\it Nature Materials}, vol.~10, November 2011.}

\subsection{Poisson's spot}

Poisson's contribution to optics
was not a successful treatment of the general phenomena of light but a prediction based on his ability to compute. Following Laplace,
he held that all phenomena could be explained by molecular interaction and was opposed to the theory of Augustin Fresnel
(1788--1827), based on a wave theory. When, in 1817, submissions for a grand prize for a study of the diffraction of
light were being examined by the Academy of Sciences, Poisson was a member of the committee in charge of examining
Fresnel's submission. Convinced that Fresnel was wrong, Poisson suggested an experiment that would prove that a
mathematical consequence of Fresnel's formulas was contrary to intuition and would disprove his theory.
When the consequence of Fresnel's theory that Poisson had derived and considered absurd was tested experimentally,
what he had judged to be absurd was actually observed: a bright spot appeared in the centre of the shadow of a disk
lit by a source situated on its axis. This phenomenon was then called derisively ``Poisson's spot''. The experiment
suggested by Poisson yielded a result in Fresnel's favor who was awarded the prize.

\subsection{The Poisson distribution}

This is probably the best known occurrence of Poisson's name in the scientific literature. In French it is referred to as
``la loi de Poisson'' (Poisson's law).

Poisson was not the first to deal with probabilities. Blaise Pascal (1623--1662), Christiaan Huygens (1629--1695), John Arbuthnot (1667--1735),
Giovanni Battista Vico (1668--1744), Georges-Louis Leclerc de Buffon (1707--1788) in his {\it Essay on Moral Arithmetic},\footnote{Essai d'arithm\'etique morale.}
of 1777, had all written on this subject, and even Voltaire had written a pamphlet, {\it
Essay on Probabilities Applied to Justice},\footnote{Essai sur les probabilit\'es en fait de justice.} in~1772, but it was not mathematical, and the revised and augmented
fifth edition of Laplace's {\it Philosophical Essay on Probabilities\footnote{Essai philosophique sur les probabilit\'es.}} of 1814
was published in volume~7 of his works in~1825.
In 1981, Bernard Bru, in his essay on Poisson and probability theory~\cite{bru}, wrote that a precursor
of the probability law that bears the name of Poisson can be found as early as 1718 in the work of the huguenot mathematician working in England,
Abraham de Moivre (1667--1754), \textit{The Doctrine of Chances}.
But in its present form, the Poisson distribution appeared for the first time on page 262 of Poisson's essay of 1829,
``Essay on the proportion of new-born females and males'',\footnote{M\'emoire sur la proportion des naissances des filles et des gar\c cons.}
which was published in the \textit{M\'emoires de l'Acad\'emie des Sciences}~\cite{naissances} of 1830, and reappeared in his subsequent book, \textit{Recherches
sur la probabilit\'e des jugements}~\cite{22}, published in~1837,
where one reads on p.~206,
\[
P= \left(1 + \omega + \frac{\omega^2}{1\cdot 2} + \frac{\omega^3}{1\cdot 2\cdot 3} + \cdots + \frac{\omega^n}{1\cdot 2\cdot 3\cdots n} \right) {\rm e}^{- \omega}.
\]

Related to the Poisson distribution are the so-called Poisson--Charlier polynomials, whose sequence first appeared in the work of the Swedish
astronomer and statistician, Carl Vilhelm Ludwig Charlier (1862--1934). It was the German mathematician
Gustav Doetsch (1892--1977) who called them Charlier polynomials in the title of an article published in \textit{Mathematische Annalen} in 1933
where he discussed the differential-difference equation they satisfied. The reviewer of Doetsch's article for \textit{Zentralblatt} wrote the
definition of these polynomials in terms of a sequence of functions defined by recursion from the Poisson distribution (Poissonsche Verteilung)
and noted the orthogonality property which they satisfy with respect to the Poisson density, whence the present terminology.
In a note in the \textit{Annals of Mathematical Statistics}
in 1947, Clifford Truesdell (1919--2000) derived their properties from those of the F-functions which he introduced, and he entitled his article,
``A note on the Poisson--Charlier functions''~\cite{truesdell}. Thus Poisson's name became attached to concepts invented a hundred years after
his death.

\subsection{The Poisson summation formula}

The Poisson summation formula is so well known that it is often called, simply, ``the Poisson formula''.
Why is the formula
\[ \sum_{n=-\infty}^{n=\infty}f(n) =
\sum_{n=-\infty}^{n=\infty}{\mathcal F}f(n),
\]
where ${\mathcal F}f$ is the Fourier transform of $f$, defined by ${\mathcal F}f(\xi) = \int_{-\infty}^{\infty} {\rm e}^{-2\pi {\rm i} x \xi} f(x) \,{\rm d}x$,
attributed to Poisson?
In the hope of finding references that would lead us to his original papers, we opened a modern textbook,
 \textit{Sphere Packings, Lattices and Groups} (1999), by John Conway and
N.J.A.~Sloane. Introducing the Jacobi theta functions, they write on p.~103, ``These functions are related by a labyrinth of identities [\dots] One may
regard [them] as consequences of the general version of the Poisson
summation formula''.

How did the reference to Poisson reach Conway and Sloane
at the end of the
20th century, and what do we find in Poisson? Fortunately, they refer to p.~475 of the
classical treatise, \textit{A~Course
of Modern Analysis} (1927), by Whittaker and Watson who in turn refer to Poisson's article of 1823 in the \textit{Journal de l'\'Ecole polytechnique},
``Continuation of the essay on definite integrals and the summation of series''\footnote{Suite
du m\'emoire sur les int\'egrales d\'efinies
et sur la sommation des
s\'eries.}~\cite{17}. There, on p.~420, we read the formula
\[
\pi + 2 \pi \sum_{n=1}^\infty {\rm e}^{-4k \pi^2n^2} = \frac{\sqrt{\pi}}{2
 \sqrt{k}} +
\frac{\sqrt{\pi}}{ \sqrt{k}} \sum_{n=1}^\infty {\rm e}^{-\frac{n^2}{4k}},
\]
which Poisson derived when working on a precise evaluation of the remainder in the summation formula that Leonhard Euler had obtained in 1736 and which Colin Maclaurin stated in his 1742 ``Treatise of Fluxions".
Since the Fourier transform of $f(x)= {\rm e}^{- \alpha x^2}$ is
${\mathcal F}f(\xi) = \sqrt{\frac{\pi}{\alpha}} {\rm e}^{-\frac{\pi^2
 \xi^2}{\alpha}}$, when the preceding formula is rewritten as
\[
\sum_{n=-\infty}^{n=\infty} {\rm
 e}^{-4 k \pi^2 n^2} =
\frac{1}{\sqrt{4 \pi k}} \sum_{n=-\infty}^{n=\infty}
{\rm e}^{-\frac{n^2}{4k}} ,
\]
we recognize instances of the summation formula,
\[
\sum_{n=-\infty}^{n=\infty}f(n) =
\sum_{n=-\infty}^{n=\infty}{\mathcal F}f(n),
\] for each function $f_k$ defined by $f_k(x)= {\rm e}^{-4k\pi^2 x^2}$.

What Whittaker and Watson observed was that, setting $4k\pi = - {\rm i}\tau$, Poisson's formula can be rewritten as
\[
\theta_3(0|\tau) = \frac{1}{\sqrt{- {\rm i} \tau}} \theta_3\left(0\,\bigg|\, {-} \frac{1}{\tau}\right),
\]
the particular case for $z =0$ of the general transformation formula for the third theta function,
\[
\theta_3(z|\tau)= (-{\rm i} \tau)^{-\frac{1}{2}} {\rm e}^{\frac{z^2}{\pi {\rm i}
 \tau}} \theta_3\left(\frac{z}{\tau}\,\bigg|\,{-}\frac{1}{\tau}\right).
\]
They also stated that a more general case is to be found in Poisson's ``Essay on the numerical calculation of definite
integrals'',\footnote{M\'emoire sur le calcul num\'erique des int\'egrales d\'efinies.} published in 1827 in the {\it
M\'emoires de l'Acad\'emie des Sciences}, a year before Jacobi published ``Continuation of the notices on elliptical
functions''\footnote{Suite des notices sur les fonctions elliptiques.} in \textit{Crelle's Journal}, which was followed by his
comprehensive treatment of the identities satisfied by the theta functions in his {\it Fundamenta nova theorias
functionum ellipticarum}, published in K\oe nigsberg in~1829.

Today, the summation formula, generalized in the theory of group representations, has applications to network theory and
error-correcting codes, that Poisson could not have anticipated.

\subsection{The Poisson kernel and the Poisson integral formula}

If one opens any modern book on potential theory, one will no doubt find a definition of ``the Poisson kernel'' and a proof of
``the Poisson integral formula'', often simply called ``the Poisson formula'', for the case of a half-plane and for a disk in the
plane, often also for the sphere in 3-space or in higher-dimensional spaces. How did these formulas reach the modern authors
and where did they appear in Poisson's vast mathematical production?

What became known as ``the Dirichlet problem'' for a domain in $n$-space, which consists of determining the value of a harmonic
function in the interior of the domain, given the value of the function on the boundary of the domain, was formulated for a~disk in the plane by Peter Gustav Lejeune Dirichlet (1805--1859) in an article in \textit{Crelle's Journal} in 1828.

In a paper on the later history of Fourier series, Jean-Pierre Kahane (1926--2017) listed, among the
advances made in the 19th century on the subject of the trigonometric series, the solution by Hermann Amandus Schwarz (1843--1921) of
the Dirichlet problem for the circle by means of the Poisson formula, in 1872.
In fact, Schwarz's article \cite{schwarz}
in the \textit{Journal f\"ur die reine und angewandte Mathematik} of 1872, contains the formula that expresses the value of a harmonic function
in the interior of a disk as an integral involving only the values of that function on the boundary circle,
\[
u(r,\phi) = \frac{1}{2\pi} \int_0^{2\pi} u(1,\phi) \frac{1-r^2}{1-2r \cos (\psi -\phi) +r^2}\, {\rm d}\psi,
\]
but he writes: ``It is easy to recognize the fundamental idea of Poisson's proof in the proof that is to be found
in Section~5.{\it b}''.\footnote{Man wird auch in dem Beweise, der in Section~5 unter~{\it b} enthalten ist, die Grundgedanken des
Poissonschen Beweises leicht wiedererkennen.} Schwarz attributed this formula to Carl Neumann (1832--1925) in his article
\cite{neumann} in volume 59 of the same journal. One does find this formula on p.~365 of Neumann's article. While
there is no mention of Poisson in Neumann's article, Schwarz on the other hand gave numerous references to Poisson's articles of 1815 and 1823,
and to his book on the theory of heat, {\it Th\'eorie math\'ematique de la chaleur}, of 1835, as well as to three other papers published in 1827,
1829 and 1831.
He thus gave us a useful map to enter the maze of sometimes very long essays that Poisson wrote, many of them with a ``Suite'' and an ``Addition''.

Poisson's search started as early as 1813 in his ``Essay on definite integrals''\footnote{M\'emoire sur les int\'egrales d\'efinies.}
\cite{15} published in the {\it Journal de l'\'Ecole polytechnique}. This paper was followed by another essay, sixty pages long, in
1820, ``Essay on the manner of expressing functions by series of periodic quantities and on the use of this transformation in the solution
of various problems''\footnote{M\'emoire sur la mani\`ere d'exprimer les fonctions par des s\'eries de quantit\'es p\'eriodiques, et sur l'usage
de cette transformation dans la r\'esolution de diff\'erents probl\`emes.}~\cite{16}, published in the same journal, and a ``Continuation
of the essay on definite integrals and the summation of series''\footnote{Suite du m\'emoire sur les int\'egrales d\'efinies et sur la sommation des
s\'eries.} \cite{17} in 1823.

Poisson was trying to establish interpolation formulas \`a la Lagrange.
In 1820, in his essay on series of periodic functions, he proved that, given a finite sequence of $m-1$ quantities, $y_1,\dots,y_{m-1}$,
setting $z_j =\frac{2}{m}\sum_{k=1}^{m-1}y_k \sin \frac{kj\pi}{m}$, it follows that $y_n=
\sum_{j=1}^{m-1} z_j \sin \frac{nj\pi}{m}.$
He wrote that this formula is a particular case of
 Lagrange's formula in his ``Researches on the nature of sound and its propagation''\footnote{Recherches sur la nature, et la propagation du
 son.} which had appeared in the first volume of the \textit{M\'emoires de l'Acad\'emie de Turin} in 1759.
Then, using a limiting procedure and exchanging summation and integration, he derived the following formula
\[
fx = \frac{2}{l}\int_0^l \sum_{k=1}^{\infty} \sin \frac{k\pi x}{m}\sin \frac{k\pi \alpha}{m} f\alpha \, {\rm d}\alpha,
\]
which he also attributed to Lagrange.
In this long essay, he expressed his aim of replacing the summation of a series by the computation of an integral or conversely,
and he treated the question of the summation of series of sines and cosines. He wrote: ``It will be advantageous to bring all of them
together under the same point of view and to deduce these values by a uniform method''.\footnote{Il ne sera pas inutile de les r\'eunir toutes
sous un m\^eme point de vue et de d\'eduire ces valeurs d'une m\'ethode uniforme.} On p.~422, he introduced the evaluation of a function
with the help of integration by way of what would later be called a Poisson kernel,
and he gave numerous applications for it, including to the motion of a vibrating string composed of two parts of different
material and the motion of a heavy body suspended from an elastic wire.

What is clear from reading Poisson is that he was not trying to solve the so-called Dirichlet problem, except maybe in the case of some
applied problem related to his more general researches, but that he returned many times to the theory of trigonometric series and that he
was in fact trying to prove the so-called Fourier theorem,
that is, he was working towards a proof of the convergence of the ``Fourier series'' of a given function to the
function itself. His attempts at a rigorous treatment of this question, as well as the later treatment by Cauchy in 1823, were not
successful, but it is in the course of such a research that Poisson introduced the function that became known as the Poisson kernel
and the integral that became known as the Poisson integral. Both terms are justified, but their appearance in the theory of harmonic
functions and potential theory came later, with Dirichlet and Schwarz.
In conclusion, we can affirm that, on the one hand, the Poisson kernel was indeed introduced by Poisson in his attempts to prove the convergence
of the Fourier series of general types of functions, and that, on the other hand, Poisson did not use the corresponding integral formula in the
search for a general solution of the Laplace equation, but only in particular cases that arose from physical problems.

\subsection{Poisson brackets}

It all began in celestial mechanics. Following Lagrange's essays of 1808 and 1809 on the variation of the principal axes of the orbits of
planets and on a general theory of the variation of arbitrary constants in mechanics, but it was in Poisson's famous essay of 1809,
``Essay on the variation of arbitrary constants in questions of mechanics''\footnote{M\'emoire sur la variation des constantes arbitraires
dans les questions de m\'ecanique.}~\cite{13} that the Poisson brackets appeared in their own right. Lagrange had derived their
expression by an involved procedure which amounted to~-- in modern terms~-- inverting the matrix of components of the canonical
symplectic 2-form. Poisson denoted the coordinates of the position of the body by $\phi$, $\psi$, $\theta$ and the components of its
velocity by $s$, $u$, $v$, and he wrote:
\begin{quote}
 It is clear that the left hand side of this equation is a complete differential with respect to $t$; integrating, we thus obtain the very
 simple equation
 \[
 {{\rm d}b \over {\rm d}s}\cdot{{\rm d}a \over {\rm d}\phi} - {{\rm d}a \over {\rm d}s}\cdot{{\rm d}b \over {\rm d}\phi}
 + {{\rm d}b \over {\rm d}u}\cdot{{\rm d}a \over {\rm d}\psi} - {{\rm d}a \over {\rm d}u}\cdot{{\rm d}b \over {\rm d}\psi}
 + {{\rm d}b \over {\rm d}v}\cdot{{\rm d}a \over {\rm d}\theta} - {{\rm d}a \over {\rm d}v}\cdot{{\rm d}b \over {\rm d}\theta}
 = \text{const}.
\]
 One sees that the constant on the right hand side of this equation will, in general, be a function of $a$ and $b$, and of the arbitrary
 constants appearing in the other integrals of the motion;
[\dots] but in order to recall the origin of this quantity, which represents a certain combination of the partial differentials of the values
of $a$ and $b$, we shall make use of the notation $(b,a)$ to denote it; so that we shall have generally\footnote{Il est visible que le premier membre de cette \'equation est une
 diff\'erentielle compl\`ete par rapport \`a $t$; en int\'egrant, nous aurons donc
 cette \'equation fort simple
\[
 {{\rm d}b \over {\rm d}s}\cdot{{\rm d}a \over {\rm d}\phi} - {{\rm d}a \over {\rm d}s}\cdot{{\rm d}b \over {\rm d}\phi}
 + {{\rm d}b \over {\rm d}u}\cdot{{\rm d}a \over {\rm d}\psi} - {{\rm d}a \over {\rm d}u}\cdot{{\rm d}b \over {\rm d}\psi}
 + {{\rm d}b \over {\rm d}v}\cdot{{\rm d}a \over {\rm d}\theta} - {{\rm d}a \over {\rm d}v}\cdot{{\rm d}b \over {\rm d}\theta}
 = \mbox{const}.
\]
 On con\c{c}oit que la constante qui fait le second membre de cette \'equation, sera
 en g\'en\'eral une fonction de~$a$ et~$b$, et des constantes arbitraires contenues
 dans les autres int\'egrales des \'equations du mouvement;
[\dots] mais, afin de rappeler l'origine de cette quantit\'e,
 qui repr\'esente une certaine combinaison des diff\'erences partielles des valeurs
 de $a$ et $b$, nous ferons usage de cette notation $(b,a)$, pour la d\'esigner;
 de mani\`ere que nous aurons g\'en\'eralement
 \[
 {{\rm d}b \over {\rm d}s}\cdot{{\rm d}a \over {\rm d}\phi} - {{\rm d}a \over {\rm d}s}\cdot{{\rm d}b \over {\rm d}\phi}
 + {{\rm d}b \over {\rm d}u}\cdot{{\rm d}a \over {\rm d}\psi} - {{\rm d}a \over {\rm d}u}\cdot{{\rm d}b \over {\rm d}\psi}
 + {{\rm d}b \over {\rm d}v}\cdot{{\rm d}a \over {\rm d}\theta} - {{\rm d}a \over {\rm d}v}\cdot{{\rm d}b \over {\rm d}\theta}
 = (b,a).
 \]}
\[
 {{\rm d}b \over {\rm d}s}\cdot{{\rm d}a \over {\rm d}\phi} - {{\rm d}a \over {\rm d}s}\cdot{{\rm d}b \over {\rm d}\phi}
 + {{\rm d}b \over {\rm d}u}\cdot{{\rm d}a \over {\rm d}\psi} - {{\rm d}a \over {\rm d}u}\cdot{{\rm d}b \over {\rm d}\psi}
 + {{\rm d}b \over {\rm d}v}\cdot{{\rm d}a \over {\rm d}\theta} - {{\rm d}a \over {\rm d}v}\cdot{{\rm d}b \over {\rm d}\theta}
 = (b,a).
 \]
 \end{quote}
(There is a short but detailed discussion of Poisson's paper in Ren\'e Dugas's {\it Histoire de la M\'ecanique}~\cite{histoiremeca}, in a
somewhat modernized notation, which renders Poisson's argument easy to understand and facilitates the reading of his original text.)

It is in this early article that Poisson introduced the change of variables from $(q_i, \dot q_i)$ to $(q_i, p_i)$, where the $p_i$'s
are the conjugate quantities, or momenta, of the $q_i$'s, defined by $\frac{\partial L}{\partial \dot q_i} =p_i$, when~$L$ is the Lagrangian
function, paving the way for the Hamiltonian form of the equations of motion, already implicit in Lagrange, eventually derived by Cauchy in a
lithographed memoir in 1831, which was only later printed, in 1834 in Italian and in 1835 in French, and eventually published by William Rowan
Hamilton (1805--1865) in his ``Second essay on a general method in dynamics'' in 1835. Poisson returned to the subject in 1816.

Jacobi, in 1850, read a biography of Poisson (maybe Arago's?), and he re-discovered what became known as Poisson's theorem, that the Poisson
bracket of two integrals of the motion is an integral of the motion. He exclaimed that this theorem was a ``truly prodigious theorem'' (ce
th\'eor\`eme vraiment prodigieux), and he endeavored to explain what he claimed its author and later authors had not perceived.

The ubiquity of the concepts of Poisson brackets, Poisson algebras and Poisson manifolds in mechanics, theoretical physics and an impressive
number of fields of mathematics suggests the question: how did this happen?
The story of Poisson brackets involves, in addition to
Lagrange, Cauchy and Hamilton, mainly Jacobi, Liouville and Sophus Lie (1842--1899), and culminates in the explanation of
the role they played in the development of quantum mechanics. It is of course too long a story to outline here.
Even a short history of Poisson geometry would imply an excursus into the history of symplectic geometry.
I shall only comment here on the fact that the name ``Poisson brackets'' does not seem to have been adopted until Whittaker used it in his
\textit{History of the Theories of {\AE}ther and Electricity} in 1910. They had previously been called ``expressions'' generally, while one author, Joseph Graindorge (1843--1889), wrote in 1872 that he was using ``Poisson's notation" (la notation de Poisson).

Already in 1857 Arthur Cayley (1821--1895), in his {\it Report of the British Association for the Advancement of Science}, had predicted the
importance that the Poisson brackets would assume later, as opposed to the Lagrange parentheses, ``The theory of Poisson gives rise to
developments which seem to have nothing corresponding to them in the theory of Lagrange''. In fact, while the Lagrange parentheses are the
components of a closed 2-form, a concept that appeared only with \'Elie Cartan (1869--1951) just before 1900, the Poisson brackets satisfy
the identity that Jacobi proved ca.\ 1840 (published posthumously in 1862), and the Jacobi identity has become the foundation of the theory of Lie
algebras.

\section{Conclusion}

Most, but not all, of the judgments that were passed on Poisson's {\oe}uvre in the 19th century were extremely laudatory.
The physicist and astronomer Fran\c cois Arago (1786--1853), who was then Secr\'etaire perp\'etuel (president) de l'Acad\'emie des Sciences,
declared at Poisson's funeral in 1840: ``Genius does not die thus; it survives in its works'',\footnote{Le g\'enie ne meurt pas ainsi; il se survit
dans ses {\oe}uvres.} and ten years later he wrote in another eulogy in the form of a scientific biography that Poisson had ``three qualities:
genius, devotion to work and mathematical erudition''.\footnote{On se demandera sans doute comment, durant une vie si courte et consacr\'ee en
grande partie au professorat, notre confr\`ere \'etait parvenu \`a attaquer et \`a r\'esoudre tant de probl\`emes. Je r\'epondrai que c'est par
la r\'eunion de trois qualit\'es: le g\'enie, l'amour du travail et l'\'erudition math\'ematique.} But the author of a history of the
mathematical and physical sciences in 12 volumes published in the 1880's,
Maximilien Marie, devoted several pages to a derogatory summary of Poisson's activity, ``Poisson
 was very far from accomplishing the promise of his youth'',\footnote{Poisson n'a pas tenu, \`a beaucoup pr\`es, les promesses de sa jeunesse.}
 claiming that his position in any debate would always have been the wrong one, and
 Pierre Costabel in~\cite{costabel}, although less negative, wrote a harsh appreciation.

When the first new buildings of the Sorbonne were inaugurated in 1890, fifty years after Poisson's death, it was Charles Hermite (1822--1901)~--
after whom the ``Hermite polynomials'' and ``Hermitian matrices'' are named~-- who delivered a speech~\cite{hermite} at the ceremony in the
presence of the then French President, Sadi Carnot.
He chose to review the work and the legacy of the professors who held the first chairs at the Facult\'e des Sciences
at its creation in 1809. He declared, ``Poisson is one of the great geometers of this century'',\footnote{Poisson est l'un des grands g\'eom\`etres
de ce si\`ecle.}
before reviewing Poisson's achievements in mathematical physics: ``For Laplace and Poisson, pure analysis is not the object, but
the tool'',\footnote{Pour Laplace et Poisson, l'Analyse pure n'est point le but, mais l'instrument.} and he pursues:
``But, having a different object, Poisson and Fourier contributed to the development of analysis which they enriched with methods, new
results, fundamental concepts''.\footnote{Mais en ayant un autre but, Poisson et Fourier contribuent au d\'eveloppement de l'Analyse
qu'ils enrichissent de m\'ethodes, de r\'esultats nouveaux, de notions fondamentales.}
 He underlined the fact that Poisson was a disciple of Laplace, but he also announced the following, referring to Poisson's work on what was
 already known as ``the Poisson theorem'',
``But he is also related to contemporary analysis regarding a question of the utmost importance and of particular interest from the point of view
of mathematical invention''.\footnote{Mais il se rattache aussi \`a l'analyse de notre temps dans une question de la plus grande importance
et qui pr\'esente un int\'er\^et singulier au point de vue de l'invention math\'ematique.} And he goes on to quote a sentence in Latin written
by Jacobi which I understand to mean, ``We here have an example that clearly shows that, if problems are not already formulated in our minds, it
may well be that we would not see most important discoveries that are set in front of our eyes''.

Was Poisson a mathematician or a physicist?
If he was called ``a geometer'' by Arago, Hermite and others, it is because ``geometer'' was the generic term for mathematician until much later in
the 19th century. His ambition was to write a comprehensive treatise of mathematical physics.
Lagrange, Laplace, Gauss, Cauchy, all contributed both to ``pure mathematics'' and to the solution of problems in physics, sometimes very practical
problems, such as geodesy, working towards the latter with the mathematial tools that they forged.
Poisson's work on the theory of magnetism had important applications to the navigation instruments for ships. His molecular theory, following
Laplace, did not win against the wave theory advocated by Fresnel, but in the history of elasticity some of his insights which had been
discredited have regained their importance in the latest development of composite materials.

Poisson never traveled. After he came to Paris at the age of 17 to enter the \'Ecole Polytechnique, he left the capital only to visit Laplace at
Arcueil, a few miles from the southern edge of Paris, where he joined the other members of the Soci\'et\'e d'Arcueil, a small circle of renowned
scientists, named after their meeting place.
But his publications abroad were numerous and his influence in Germany, in England, and in Russia was considerable.
Several of his books were translated and extracts or summaries of his articles appeared in the \textit{Zeitschrift f\"ur Physik},
in the {\it Annalen der Physik}, and the {\it Philosophical Transactions}.\footnote{These references
are not in the autograph list of his works~-- now in the Biblioth\`eque de l'Institut~-- that Poisson himself drew up not long before he died.
They had escaped the meticulous, invaluable work of Pierre Dugac in~\cite{dugac}.}

But Lie algebra theory had to wait for other geniuses to be developed, and Poisson geometry had to wait for more than a century and a half to be
developed in various forms in the work of, among others, Wolfgang Pauli, George W. Mackey, Wlodzimierz Tulczyjew, Vladimir Maslov, Robert Hermann,
Alexander Kirillov, Mosh\'e Flato, Andr\'e Lichnerowicz and Alan Weinstein.

Poisson's role in French science was dominant while he lived.
His explanations of physical phenomena were mostly proved wrong by his contemporaries or by later physicists, but his achievements in mathematical
physics remain. Often, he found the right equations using the wrong physical arguments.
He was a formidable ``computer'' and his legacy in mathematics is essential.
He advanced mathematics by trying to solve physical problems, sometimes successfully, sometimes not, and he is to be remembered for concepts,
equations, formulas and theorems. He demonstrated how physical problems can suggest entirely abstract mathematics, what we now call mathematical
physics.

\pdfbookmark[1]{References}{ref}
\LastPageEnding

\end{document}